\documentclass[11pt]{article}

\usepackage{amsmath,amssymb}
\usepackage{graphicx}  
\usepackage{pict2e}
\newcommand{\Ex}{\mathbb{E}}
 \renewcommand{\Pr}{{\mathbb{P}}}
 \renewcommand{\SS}{\mathcal S}
 \newcommand{\LL}{\mathcal L}
\newcommand{\GG}{\mathcal G}
\newcommand{\Ints}{{\mathbb{Z}}}
 
\newtheorem{Lemma}{Lemma}

\newtheorem{Proposition}[Lemma]{Proposition}
\newtheorem{Corollary}[Lemma]{Corollary}
\newtheorem{Conjecture}[Lemma]{Conjecture}
\newtheorem{OP}[Lemma]{Open Problem}
 \newcommand{\proof}{{\bf Proof.\ }}
 \newcommand{\qed}{\ \ \rule{1ex}{1ex}}
 \newcommand{\var}{\mathrm{var}}
\begin{document}

\title{The Nearest Unvisited Vertex Walk on Random Graphs}
 \author{David J. Aldous\thanks{Department of Statistics,
 367 Evans Hall \#\  3860,
 U.C. Berkeley CA 94720;  aldous@stat.berkeley.edu;
  www.stat.berkeley.edu/users/aldous.}
}

 \maketitle
 
 \begin{abstract}
 We revisit an old topic in algorithms, the deterministic walk on a finite graph which always moves toward the nearest unvisited vertex until every vertex is visited.
 There is an elementary connection between this cover time and ball-covering (metric entropy) measures. 
 For some familiar models of random graphs,  this connection
 allows the order of magnitude of the cover time to be deduced from first passage percolation estimates. 
 Establishing sharper results seems a challenging problem.
 \end{abstract}
 

\section{Introduction}
Consider a connected undirected graph $G$ on $n$ vertices, where the edges $e$ have positive real lengths $\ell(e)$.
Consider an entity -- let's call it a robot -- that can move at speed $1$ along edges.
There are many different rules one might specify for how the robot chooses which edge to take after reaching a vertex  
-- for instance the ``random walk" rule, to choose edge $e$ with probability proportional to $\ell(e)$ or $1/\ell(e)$. 
One well-studied aspect of the random walk is the {\em cover time}, the time until every vertex has been visited -- see 
Ding, Lee and Peres \cite{cover} for references to 
special examples and surprisingly deep connections with other fields.
This article instead concerns what we will call\footnote{Confusingly previously called {\em nearest neighbor}, inconsistent with the usual terminology that neighbors are linked by a single edge, but justifiable by the artifice of extending the given graph to a complete graph via defining each edge $(v,v^*)$ to have length $d(v,v^*)$.  But the phrase {\em nearest neighbor} is used in many other contexts, so the more precise name NUV seems preferable.}
the {\em nearest unvisited vertex} (NUV) walk, defined as follows.
A path of edges has a length, the sum of edge-lengths, and the distance $d(v,v^*)$ 
between vertices is the length of the shortest path.
For simplicity assume all such distances are distinct, so the shortest path is unique.
Now the NUV walk is the deterministic walk defined in words by
\begin{quote}
after arriving at a vertex, next move at speed $1$ along the path to the closest unvisited vertex
\end{quote}
and continue until every vertex has been visited.\footnote{This {\em walk} convention is consistent with random walk cover times; one could alternatively
use the {\em tour} convention that the walk finally returns to its start, consistent with TSP.}
In symbols, from initial vertex $v_0$ 
the vertices can be written $v_0,v_1,v_2, \ldots,v_{n-1}$ in order of first visit;
\begin{equation}
 v_i = \arg \min_{v \not\in \{v_0,\ldots,v_{i-1}\}} d(v_{i-1},v) , \quad 1 \le i \le n-1
 \label{vii}
 \end{equation}
and this walk has length $L = L_{NUV}  = L_{NUV}(G,v_0) = \sum_{i=1}^{n-1} d(v_{i-1},v_i)$.

There are several types of question one can ask about NUV walks.
\begin{itemize}
\item The order of magnitude of $L$ for a general graph?
\item Sharper estimates of $L$  for specific models of random graphs?
\item Structural properties of the NUV path in different contexts?
\end{itemize}
The first question has  been studied in the context of  TSP (travelling salesman problem) heuristics and  robot motion, 
and a 2012 survey of the general area, under the name {\em online graph exploration}, is given in 
Megow, Mehlhorn and Schweitzer \cite{megow}.

\subsection{Outline of results}
Our first purpose is to record a formalization (Proposition \ref{P:1}) of the basic general relationship between $L_{NUV}$ and ball-covering.
This is  implicit in two now-classical results:
Corollary \ref{C:1}, which compares $L_{NUV}$ to the length $L_{TSP}$ of the shortest path through all $n$ vertices,
and Corollary \ref{C:2}, which upper bounds $L_{NUV}$ for $n$ arbitrary points in the unit square with Euclidean distance.
As shown in section \ref{sec:basic}, each follows easily from our formalization.

Our main  purpose is to point out that the relation with ball-covering enables (in some simple probability models)
the order of magnitude of $L$ to be deduced easily from known  first passage percolation estimates.
In section \ref{sec:FPP} we study two specific models.
\begin{itemize}
\item 
For the $m \times m$ grid with i.i.d. edge-lengths, Corollary \ref{C:grid} shows that $L$ is indeed $O(m^2)$ rather than larger order.
\item
For the complete graph on $n$ vertices, with i.i.d. edge-lengths normalized so that the shortest edge at a vertex is order $1$, 
Corollary \ref{C:MF} shows that $L$ is indeed $O(n)$ rather than larger order. 
\end{itemize}
In both of those models the (first-order) behavior of  first passage percolation is well understood, via the {\em shape theorem} on the two-dimensional grid, 
and the Yule process approximation on the complete graph model.

A final purpose is to point out that the second and third questions above have apparently never been studied.
The NUV rule on a deterministic graph is ``fragile" in the sense that small changes in the length of an edge might affect a large proportion of the walk,
But it is possible that introducing random edge-lengths might ``smooth" the typical properties of the walk on a random graph. 
We defer further general discussion to section \ref{sec:remarks}.

\section{Basics}
\label{sec:basic}
\subsection{Relation with ball-covering}
A basic mathematical observation is that $L_{NUV}$ is  related to ball-covering\footnote{And thereby to {\em metric entropy} -- see section \ref{sec:order}}.
Given $r>0$ define $N(r) = N(G,r)$ to be the minimal size of a set $\SS$ of vertices such that every vertex is within distance $r$ from 
some element of $\SS$.  
In other words,  the union over $s \in \SS$ of the balls of radii $r$ centered at $s$ covers the entire graph.
\begin{Proposition}
\label{P:1}
(i) $N(r) \le 1 + L_{NUV}/r, \ 0 < r < \infty $.\\
(ii) $L_{NUV}  \le 2 \int_0^{\Delta/2} N(r) \ dr $ where 
 $\Delta = \max_{v,w} d(v,w)$ is the diameter of the graph.
 \end{Proposition}
 \proof
Inequality (i) is almost obvious.
As at (\ref{vii}), write the vertices as $v_0,v_1,v_2, \ldots,v_{n-1}$ in order of first visit by the NUV walk, and say $v_i$ has rank $i$.
Write $\zeta(v_i) = \sum_{j=0}^{i-1} d( v_j, v_{j+1})$ for the length of the walk up to $v_i$.
Select vertices $(z(k), 0 \le k \le k^* - 1)$ along the walk by selecting the first vertex at distance $>r$ along the walk after the previous selected vertex.
That is, $z(k) = v_{I(k)}$ where $I(0) = 0$ and for $k \ge 0$ 
\[ I(k +1) = \min \{i > I(k) : \zeta(v_i)  - \zeta(v_{I(k)} ) > r \}  \]
until no such $i$ exists.
By construction every vertex is within distance $r$ of some $z$, and the number $k^*$ of selected vertices is at  most
$1 + L_{NUV}/r$.
This establishes (i).

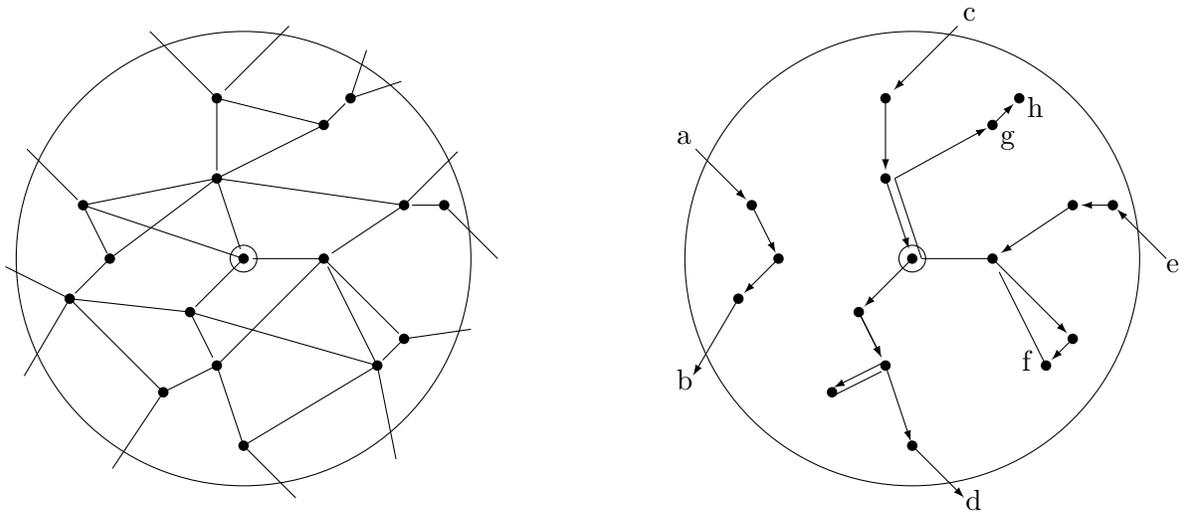
\begin{figure}

\setlength{\unitlength}{0.14in}
\begin{picture}(20,20)(-5,-10)
\put(0,0){\circle{17}}
\put(0,0){\circle{1}}

\put(-8.1,4.1){\line(1,-1){1.9}}
\put(-6,2){\line(1,-2){0.9}}
\put(-5,0){\line(-1,-1){1.3}}
\put(-6.5,-1.5){\line(-1,-1.7){1.7}}

\put(-6,2){\circle*{0.39}}
\put(-5,0){\circle*{0.39}}
\put(-6.5,-1.5){\circle*{0.39}}

\put(1.7,8.7){\line(-1,-1){2.4}}
\put(-1,6){\line(0,-1){2.7}}
\put(-1,3){\line(1,-3){0.88}}
\put(0,0){\line(-1,-1){1.8}}
\put(-2,-2){\line(1,-2){0.85}}
\put(-1,-4){\line(-2,-1){1.85}}
\put(-1,-4){\line(1,-3){0.95}}
\put(0,-7){\line(1,-1){1.95}}

\put(-1,6){\circle*{0.39}}
\put(-1,3){\circle*{0.39}}
\put(0,0){\circle*{0.39}}
\put(-2,-2){\circle*{0.39}}
\put(-1,-4){\circle*{0.39}}
\put(-3,-5){\circle*{0.39}}
\put(0,-7){\circle*{0.39}}

\put(9.5,0){\line(-1,1){1.85}}
\put(7.5,2){\line(-1,0){1.2}}
\put(6,2){\line(-3,-2){2.7}}
\put(3,0){\line(1,-1){2.8}}
\put(6,-3){\line(-1,-1){0.8}}
\put(5,-4){\line(-1,2){1.85}}
\put(3,0){\line(-1,0){2.65}}
\put(-1,3){\line(2,1){3.85}}
\put(3,5){\line(1,1){0.8}}

\put(7.5,2){\circle*{0.39}}
\put(6,2){\circle*{0.39}}
\put(3,0){\circle*{0.39}}
\put(6,-3){\circle*{0.39}}
\put(5,-4){\circle*{0.39}}

\put(3,5){\circle*{0.39}}
\put(4,6){\circle*{0.39}}

\put(-1,6){\line(-1,1){2.5}}
\put(-1,6){\line(4,-1){4}}
\put(4,6){\line(1,3){0.6}}
\put(4,6){\line(3,1){1.9}}
\put(-1,3){\line(-5,-1){5}}
\put(-1,3){\line(-4,-3){4}}
\put(-1,3){\line(7,-1){7}}  
\put(6,2){\line(1,1){2.0}}  
\put(0,0){\line(-3,1){6}}  
\put(-2,-2){\line(-9,1){4.5}}  
\put(-6.5,-1.5){\line(-2,1){2.4}}  
\put(-6.5,-1.5){\line(1,-1){3.5}}  
\put(-3,-5){\line(-2,-3){1.9}}  
\put(-1,-4){\line(1,1){4}}  
\put(5,-4){\line(-7,2){7}}  
\put(5,-4){\line(-5,-3){5}}  
\put(5,-4){\line(1,-5){0.7}}  
\put(6,-3){\line(7,1){2.5}}

\put(25,0){\circle{17}}
\put(25,0){\circle{1}}

\put(16.9,4.1){\vector(1,-1){1.9}}
\put(19,2){\vector(1,-2){0.9}}
\put(20,0){\vector(-1,-1){1.3}}
\put(18.5,-1.5){\vector(-1,-1.7){1.7}}

\put(19,2){\circle*{0.39}}
\put(20,0){\circle*{0.39}}
\put(18.5,-1.5){\circle*{0.39}}

\put(26.7,8.7){\vector(-1,-1){2.4}}
\put(24,6){\vector(0,-1){2.7}}
\put(24,3){\vector(1,-3){0.88}}
\put(25,0){\vector(-1,-1){1.8}}
\put(23,-2){\vector(1,-2){0.85}}
\put(23.9,-3.9){\vector(-2,-1){1.85}}
\put(23,-2){\vector(1,-2){0.85}}
\put(22.1,-5.1){\line(2,1){1.75}}
\put(24,-4){\vector(1,-3){0.95}}
\put(25,-7){\vector(1,-1){1.95}}

\put(24,6){\circle*{0.39}}
\put(24,3){\circle*{0.39}}
\put(25,0){\circle*{0.39}}
\put(23,-2){\circle*{0.39}}
\put(24,-4){\circle*{0.39}}
\put(22,-5){\circle*{0.39}}
\put(25,-7){\circle*{0.39}}

\put(34.5,0){\vector(-1,1){1.85}}
\put(32.5,2){\vector(-1,0){1.2}}
\put(31,2){\vector(-3,-2){2.7}}
\put(28,0){\vector(1,-1){2.8}}
\put(31,-3){\vector(-1,-1){0.8}}
\put(30,-4){\line(-1,2){1.75}}
\put(28,0){\line(-1,0){2.65}}
\put(25.35,0){\line(-1,3){1.0}}
\put(24.35,3){\vector(3.65,2){3.45}}
\put(28,5){\vector(1,1){0.8}}

\put(32.5,2){\circle*{0.39}}
\put(31,2){\circle*{0.39}}
\put(28,0){\circle*{0.39}}
\put(31,-3){\circle*{0.39}}
\put(30,-4){\circle*{0.39}}

\put(28,5){\circle*{0.39}}
\put(29,6){\circle*{0.39}}

\put(16.2,4.3){a}
\put(16.2,-4.9){b}
\put(26.9,8.9){c}
\put(27,-9.4){d}
\put(34.5,-0.5){e}
\put(29.1,-4.2){f}
\put(28.3,4.3){g}
\put(29.3,5.3){h}

\end{picture}
\caption{Illustration of the proof of (\ref{newq}).  
The left panel shows the subgraph within a radius-$r$ ball.  
The NUV walk must consist of one or several excursions within the ball.
These excursions depend on the configuration outside the ball, and the right side shows one possibility.  
The first excursion enters via edge $a$ and exits via edge $b$.
The second excursion enters via edge $c$ and exits via edge $d$, en route backtracking across one edge.
The third excursion enters via edge $e$ and proceeds to vertex $f$; at that time only vertices $g, h$ within the ball are unvisited, and the next step of the walk is a path going via three previously-visited vertices to reach $g$ and then $h$.  The next step from $h$, not shown, might be very long, depending on whether nearby vertices outside the ball have all been visited.
Arrowheads indicate the end of a step of the NUV walk, that is the edge by which the vertex is first entered.
}
\label{Fig:1}
\end{figure}

For inequality (ii), write 
$D(v_i) = d(v_i,v_{i+1})$ for the length of the {\em path} (which may encompass several edges)
from the rank-$i$ vertex to the rank-$(i+1)$ vertex, and $D(v_{n-1}) = 0$.
The argument rests upon the following  simple observation, illustrated in  Figure \ref{Fig:1}.
Fix a vertex $v^*$ and a real $r > 0$, and consider the set of vertices within distance $r$ from $v^*$:
\[ B(v^*,r) := \{v : d(v,v^*) \le r \} . \]
Consider the vertex $\bar{v}$ of highest NUV-rank within $B(v^*,r)$.
When the NUV walk first visits $v_i \in B(v^*,r)$ with $v_i \neq  \bar{v}$, 
there is then some first unvisited vertex $\tilde{v}$ on the minimum-length path from $v_i$ to $\bar{v}$, and so 
\[ D(v_i) \le d(v_i,\tilde{v}) \le d(v_i,\bar{v}) \le 2r \]
the final inequality using the triangle inequality via $v^*$.
We conclude that
\begin{equation}
\mbox{ $D(v) \le 2r$ for all $v \in B(v^*,r)$ except perhaps one vertex}. 
\label{newq}
\end{equation}
Now by considering a set, say $S(r)$, containing $N(r)$ vertices, such that every vertex is within distance $r$ from some element of $S(r)$, 
inequality (\ref{newq}) implies
\begin{equation}
\mbox{the number of vertices $w$ with $D(w) > 2r$
is at most $N(r)$. }
\label{eq}
\end{equation}
Because $D(w)$ is bounded by the graph diameter $\Delta$, 
for a uniformly random vertex $J$ we have
\begin{eqnarray*}
 L_{NUV} &=&  n \Ex[D(J)] \\
 &= & n \int_0^\Delta P(D (J)>r)dr \\
&= & \int_0^\Delta \mbox{(number of vertices $w$ with $D(w) > r$)} \ dr \\
&\le& \int_0^\Delta \ N(r/2) dr  
\end{eqnarray*}
which is equivalent to (ii).
\qed

\paragraph{Remarks.}
The simple formulation of Proposition \ref{P:1} is more implicit than explicit in the literature we have found.  
Part (i) is a less sharp version of a more complex lemma used in 
Rosenkrantz, Stearns and Lewis \cite{rosen} to prove Corollary \ref{C:1} below.  
In the context of TSP or robot exploration heuristics, the NUV algorithm is typically (e.g. in  
Hurkens and Woeginger  \cite{hurkens} 
 and in
 Johnson and Papadimitriou \cite{johnson}) 
 mentioned only briefly before continuing to better algorithms.
From an algorithmic viewpoint, calculating $N(r)$ on a general graph is not simple, so part (ii) of Proposition \ref{P:1} is not so relevant, 
but as we see in section \ref{sec:FPP} it is very helpful in providing order-of-magnitude bounds for familiar models of random networks.

\subsection{Two classical results}
Two classical results follow readily from the formulation of Proposition \ref{P:1}.
Write $L_{TSP}  = L_{TSP}(G,v_0) $ for
the length of the shortest {\em walk} starting from $v_0$ and visiting every vertex\footnote{The convention that TSP refers to a {\em tour} has the virtue that the length is independent of starting vertex.  But the latter is not true for the NUV tour.}.
So $L_{NUV}  \ge L_{TSP}$ and it is natural to ask how large the ratio can be. 
This was answered in Rosenkrantz et al. \cite{rosen}.
\begin{Corollary}
\label{C:1}
Let $a(n)$ be the maximum, over all connected $n$-vertex graphs with edge lengths and all initial vertices, of the ratio
$L_{NUV}/L_{TSP}$.
Then $a(n) = O(\log n)$.
\end{Corollary}
\proof
The argument for Proposition \ref{P:1}(i) is unchanged if we use the TSP path instead of the NUV path,
so in fact gives the stronger result
 $N(r) \le 1 + L_{TSP}/r, \ 0 < r < \infty $.
 Now apply Proposition \ref{P:1}(ii) and note that  $\Delta \le  L_{TSP}$, so 
 \[ L_{NUV}  \le 2 \int_0^{L_{TSP}/2} \min(n,    1 + L_{TSP}/r )      \ dr  \le 2 L_{TSP} + 2  L_{TSP} \log n \]
 the second inequality by splitting the integral at $r =  L_{TSP}/n$.
\qed

There are examples to show that the $O(\log n)$ bound  cannot be improved -- see  
 Johnson and Papadimitriou \cite{johnson},
Hurkens and Woeginger \cite{hurkens},
Hougardy and Wilde \cite{hougardy},
Rosenkrantz et al. \cite{rosen}.
As noted in the elementary expository article Aldous \cite{me-ES}, 
in constructing such an example the key point is to make the bound in (\ref{newq}) be tight, in the sense
\begin{quote}
for appropriate values of $r$ with $1 \ll L_{TSP}/r \ll n$ there are distinguished vertices separated by distance $r$  along the TSP path such that 
the NUV path from one to the next is order $r$.
\end{quote}
Hurkens and Woeginger \cite{hurkens} show that one can make  such examples be planar, embedded in the plane with edge-lengths as Euclidean length, and edge-lengths constrained to a neighborhood of $1$.
But such constructions seem very artificial.

Here is the second classical result. 
See Steele \cite{steele} for one proof and the early history of this result.
\begin{Corollary}
\label{C:2}
There is a constant $A$ such that,
for the complete graph on $n$ arbitrary points in the unit square, with Euclidean lengths,
\[ L_{NUV}  \le A n^{1/2} . \]
\end{Corollary}
Note this implies the  well known corresponding result  $L_{TSP}  \le A n^{1/2}$ .

\proof
By ball-covering in the continuum unit square
there is a numerical constant $C$ such that $N(r) \le C/r^2$, and so Proposition \ref{P:1}(ii) gives
\[ L_{NUV}  \le 2 \int_0^{\sqrt{1/2}} \min (n, C/r^2) \ dr 
\le 4 C^{1/2} n^{1/2} .\]
\qed

\subsection{The order of magnitude question}
\label{sec:order}
What is the size of $L_{NUV}$  for a {\em typical} graph?
That is a very vague question, but let us attempt a discussion anyway.
For this informal discussion it is convenient to scale distances so that the typical distance from a vertex to its closest neighbor is order $1$, and therefore   $L_{NUV}$
is at least order $n$.
Examples mentioned above show that $L_{NUV}$ can still be as large as order $n \log n$, 
but intuition suggests that for natural examples 
 $L_{NUV}$  is of order $n$ rather than larger order.
For this it is certainly necessary, but not sufficient, that the length $L_{MST}$ of the minimum spanning tree (MST)\footnote{Recall $L_{MST} \le L_{TSP} \le 2 L_{MST}$.} 
is $O(n)$.
Proposition \ref{P:1}(ii) provides a quantitative criterion:
it is sufficient that $N(r)/n$ is order $r^{- \alpha}$ for some $\alpha > 1$ over $1 \ll r \ll \Delta$.
Intuitively this corresponds to ``dimension $> 1$", where dimension is measured by metric entropy\footnote{The reader may be more familiar with metric entropy involving {\em small} balls for continuous spaces, but it is equally relevant in our context of large balls, as used for instance in defining fractal dimension of subsets of $\Ints^d$.}, 
 as illustrated in the examples in section \ref{sec:FPP}.

\subsection{Other questions in the deterministic setting}
It is not clear what other results might hold for general graphs $G$.
One can ask about the variability of $L_{NUV}(G,v)$ as $v$ varies.
Clearly it can be arbitrarily concentrated  e.g. on the complete graph with edge-lengths arbitrarily close to $1$.
On the other hand, consider the linear graph $G_n$ on vertices $\{0,1,\ldots,n-1\}$ with slowly decreasing edge-lengths $\ell(i-1,i) = 1 - i/n^2$.
Here there is a factor of $2$ variability in $L_{NUV}(G,v)$ as $v$ varies.
We do not see any easy example with large variability, prompting the following question. 
\begin{OP}
\label{OP:1}
Is $ \frac { \max_v L_{NUV}(G,v)}{ \min_v L_{NUV}(G,v)}$ 
 bounded over all finite graphs $G$?
\end{OP}
In this context it is perhaps more natural to extend the NUV walk to a {\em tour} which finally returns to its start.
Note that in the  linear graph example above,  $|L_{NUV}(G,v) -  L_{NUV}(G,v^\prime)  |$ is small for adjacent vertices $(v,v^\prime)$,
so one can ask whether there there is a general bound for some average of $|L_{NUV}(G,v) -  L_{NUV}(G,v^\prime)  |$ over nearby vertex-pairs $(v,v^\prime)$.

One can also consider overlap of edges used in walks from different starts.
Note that if two vertices are each other's nearest neighbor then every NUV walk uses their linking edge.
One can ask, for the two walks started at arbitrary different vertices, how small 
can be the proportion of time spent on edges used by both walks,
though we hesitate to formulate a conjecture.

\subsection{The three levels of randomness}
\label{sec:3levels}
Introducing randomness leads to different questions.
There are three ways one can introduce randomness.
One can simply randomize the starting vertex.
This suggests the following conjecture, modifying Open Problem \ref{OP:1}.
\begin{Conjecture}
\label{Con:2}
The ratio $\frac{ {\rm s.d.}(L_{NUV}(G,V))}{\Ex L_{NUV}(G,V)}$, 
where the initial vertex $V$ is uniform random,
 is bounded over all finite graphs.
\end{Conjecture}

A second level of randomness is to start with a given deterministic $G$ but then consider the random graph $\GG$
in which the edge-lengths $\ell(e)$ are replaced by independent random lengths $\ell^*(e)$ 
with Exponential(mean $\ell(e)$) distribution. 
So here we have a random variable $\LL^*(G) = L_{NUV}(\GG,V)$ where again the initial vertex $V$ is uniform random.
In this model of random graphs $\GG$, results of Aldous \cite{me-FPP} for first passage percolation say that the percolation time is 
weakly concentrated\footnote{As in the weak law of large numbers.} 
around its mean provided no single edge contributes non-negligibly to the total time.
So one can ask whether a similar result holds for $\LL^*(G)$.

The third level of randomness involves more specific models 
of random graphs, which we will consider in the next sections.

\section{Random points in the square}
\label{sec:square}
\begin{figure}
\includegraphics[width=5.0in]{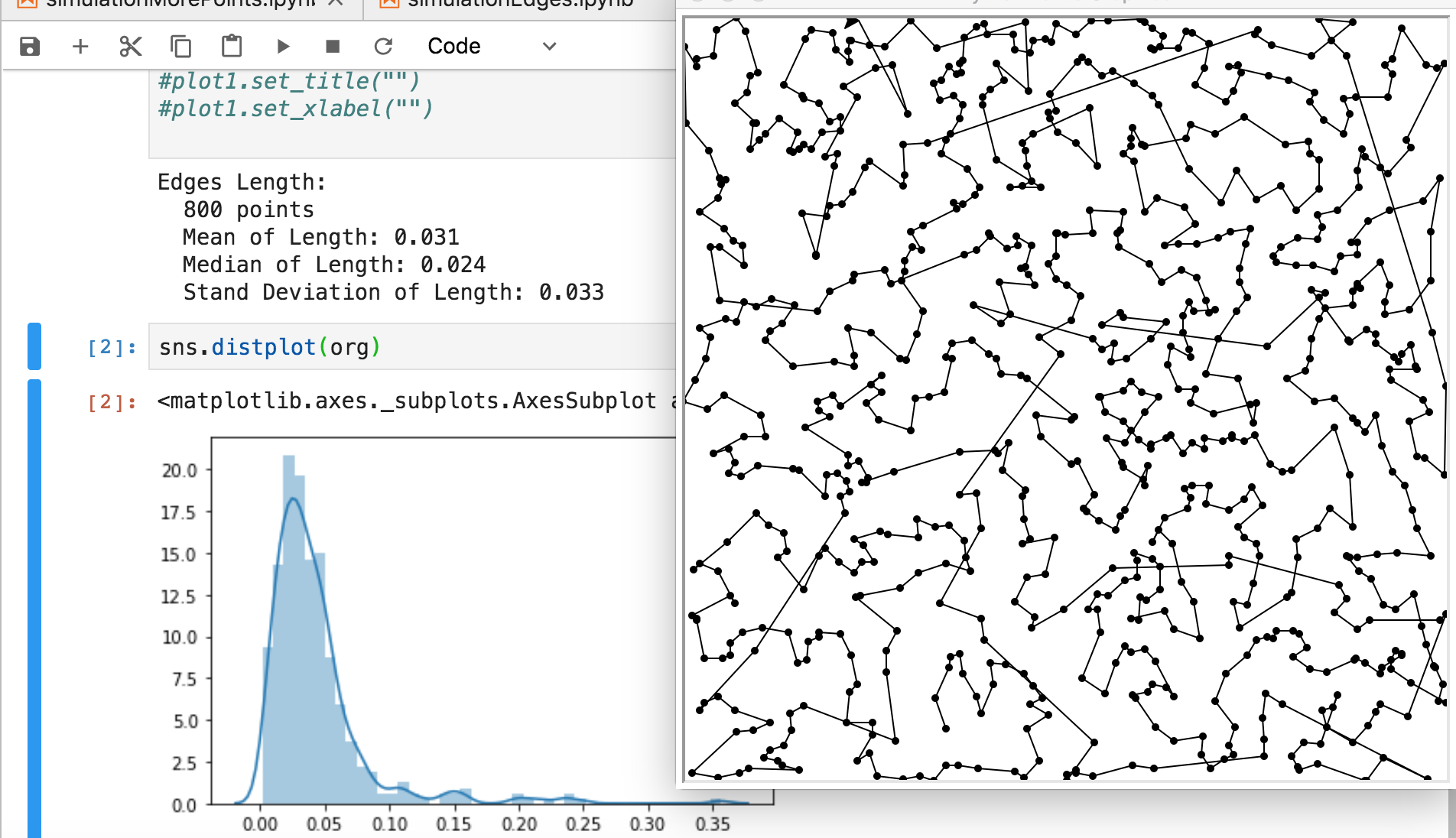}
\caption{A NUV walk through 800 random points in the unit square, and histogram of step lengths.}
\label{Fig:800}
\end{figure}

One very special model of random graph is to take the complete graph on $n$  random (i.i.d. 
uniform) points in the unit square,
with Euclidean edge-lengths.
Figure \ref{Fig:800} shows a realization of the corresponding NUV walk with $n = 800$ random points,
and Table \ref{table:3} shows some simulation data for the lengths $L^*_n$ of the NUV walk (see discussion below).
The qualitative behavior seen in simulations corresponds to intuition:
the walk starts to traverse through most (but not all) vertices in any small region, goes through different regions as some discrete analog of a space-filling curve, 
and near the end has to capture missed patches and the remaining isolated unvisited vertices via longer steps across already-explored regions.
Indeed in Figure \ref{Fig:800} we see that the actual behavior of the walk within a medium-sized ball is like the sketch in 
Figure \ref{Fig:1}, with several different excursions.

\begin{table}[h!]
\centering
\begin{tabular}{rrcc}
$n$ & $\Ex L^*_n$&$n^{-1/2} \Ex L^*_n$& s.d.($L^*_n$) \\
100 & 9.05 & 0.91 & 0.41  \\
200 &  12.78 & 0.90  & 0.54  \\
400 &  18.06 & 0.90  & 0.54 \\
800 & 25.54& 0.90 & 0.49
\end{tabular}
\caption{Simulation data for lengths $L^*_n$ in the random points in unit square model.
 Simulations and data  in this model by  Yechen Wang.}
\label{table:3}
\end{table}

The lack of scaling for the s.d. may seem surprising, but is understandable as follows.
To adhere to our scaling convention
(distance to nearest neighbor is order $1$)
we should take the square to have area $n$ and write $L_n = n^{1/2}L^*_n$ for the length of the NUV walk.
Intuition, thinking of $L_n$ as the sum of $n$ order-$1$ lengths, suggests there are limit constants
\begin{equation}
c := \lim_n n^{-1} L_n = \lim_n n^{-1/2} L^*_n; \quad \sigma := \lim_n n^{-1/2} \mathrm{s.d.} (L_n) =  \lim_n  \mathrm{s.d.} (L^*_n) .
\label{c3lim}
\end{equation}
Our small-scale simulation data suggests this holds in the present model with $c \approx 0.9$ and $\sigma \approx 0.5$.
How generally this holds is a natural question, and 
we defer further discussion to section \ref{sec:remarks}.

Corollary \ref{C:2} implies $\Ex L_n \le An$, which is all that we know rigorously.
But there are many questions one can ask.
As well as the limits (\ref{c3lim})
one might conjecture there are concentration bounds and a Gaussian limit for $n^{-1/2} (L_n  - \Ex L_n)$.
For TSP length, existence of a limit constant is known via subadditivity arguments 
(Steele \cite{steelebook} and Yukich \cite{yukich})
and concentration via now-classical Talagrand arguments, and for MST length the Gaussian limit
is also known by martingale arguments (Kesten and Lee \cite{kestenMST}).
Alas it seems hard to find any rigorous such arguments for the NUV walk.
One might also bear in mind that, for the {\em random walk} cover time problem, the two-dimensional case is the hardest to analyze sharply, so this might also hold for the NUV walk.

\begin{figure}
\begin{center}
\includegraphics[width=2.5in]{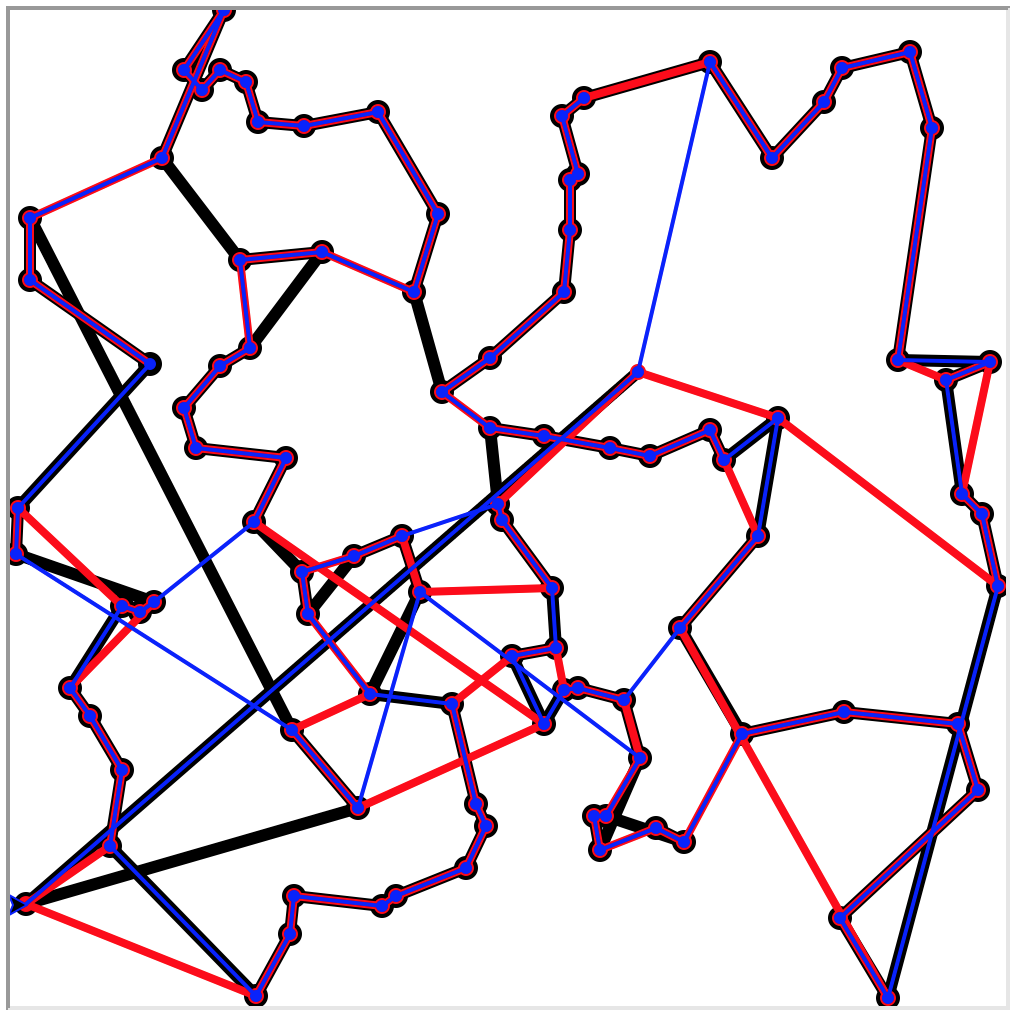}
\end{center}
\caption{3 different starts for the NUV walk on 100 points in the square.}
\label{Fig:3starts}  
\end{figure}

In any of our models, 
by considering the length as $L_n(G_n,V_n)$  for a uniform random starting vertex $V_n$,
we can consider the variance decomposition
\[
\var L_n = \var \Ex(L_n \vert G_n) + \Ex \var(L_n \vert G_n) 
\]
where the first term represents the variability due to the random graph
and  the second term represents the variability due to the starting vertex. 
In simulations of the present model, for $n = 100$ the two terms are roughly equal.  
Figure \ref{Fig:3starts}  superimposes the NUV walks from three different starts, in a realization of the present model, giving some impression of the extent of overlap.

\section{Relation with first passage percolation}
\label{sec:FPP}
For graphs with i.i.d. random edge-lengths, 
one can seek to find the correct order of magnitude of $L_{NUV}$ by
combining Proposition \ref{P:1}(ii) with known {\em first passage percolation} (FPP) results.  
Here is the basic  example.

\subsection{The 2-dimensional grid}
Consider the $m \times m$ grid, that is the subgraph of the Euclidean lattice $\Ints^2$, and assign i.i.d. edge-lengths
$\ell(e) > 0$ to make a random graph $G_m$.
Because the shortest edge-length at a given vertex is $\Omega(1)$, clearly $L_{NUV}$ is $\Omega(m^2)$.
\begin{Corollary}
\label{C:grid}
For the 2-dimensional grid model $G_m$ above, the sequence  $(m^{-2} L_{NUV}(G_m), \ m \ge 2)$ is tight.
\end{Corollary}
We conjecture that in fact $m^{-2} L_{NUV}(G_m)$ converges in probability to a constant, but we do not see any simple argument.  
Table \ref{table:7} shows simulation data, where $\ell(e)$ has Exponential(1) distribution.

\begin{table}[h!]
\centering
\begin{tabular}{rrrcc}
$n= m^2$ & $\Ex L(G_m)$&$n^{-1} \Ex L(G_m)$& s.d.($L(G_m)$) & $n^{-1/2} $ s.d.($L(G_m)$) \\
100 & 66.2 & 0.66 & 7.67 & 0.77 \\
400 &  259 & 0.65 & 14.8 & 0.74 \\
900 &  576 & 0.64 & 17.0 & 0.57
\end{tabular}
\caption{Simulation data for lengths $L(G_m) $ in the grid model.}
\label{table:7}
\end{table}

\proof
For a vertex $v$ of $G_m$ write $B(v,r)$ for the random set of vertices $v^\prime$ with $d(v,v^\prime) \le r$,
and write $D(v,r)$ for the non-random set of vertices $v^\prime$ with Euclidean distance $|| v - v^\prime|| \le r$.
Standard results for FPP on $\Ints^2$ going back to Kesten \cite{kesten} 
(see Auffinger, Damron and Hanson \cite{auff} Theorem 3.41 for recent discussion)
imply that there exist constants $c_1, c_2, c_3$ (depending on the distribution of $\ell(e)$) such that
\begin{equation}
\Pr( D(v,r) \not\subseteq B(v,c_1r)) \le c_2 \exp(- c_3r) , \ 0 < r < \infty .
\label{subseteq}
\end{equation}
The remainder of the proof is conceptually straightforward.
Given large $m$ and $r$, there is a set $S(m,r)$ of at most $a_1 m^2/r^2$ vertices of $G_m$ such that
$\cup_{v \in S(m,r)} D(v,r)$ covers $G_m$, and note $D(v,r)$ contains at most $a_2r^2$ vertices; here $a_1$ and $a_2$ are absolute constants.
By Markov's inequality and (\ref{subseteq}) the 
probability of the event
\begin{eqnarray}
 &&\mbox{the number of $v$ in $S(m,r)$ such that $D(v,r) \not\subseteq B(v,c_1r) $} \nonumber\\
&&  \mbox{ exceeds a given $s > 0$} \label{event}
\end{eqnarray}
 is at most 
$a_1 m^2 r^{-2} c_2 \exp(- c_3r) /s  $.
Apply this with $s =m^2 r^{-2}  \exp(-c_3 r/2)$.
Now define a vertex-set $S^+(m,r)$ as
\begin{quote}
 the union of $S(m,r)$ and all the vertices in all the discs $D(v,r)$ with $v \in S(m,r)$ and
$D(v,r) \not\subseteq B(v,c_1r) $.
\end{quote}
Outside the event (\ref{event}), we have that
$\cup_{v \in S^+(m,r)} D(v,r)$ covers $G_m$, and $S^+(m,r)$ has cardinality at most
\[ n_m(r) := a_1 m^2/r^2 + s a_2r^2 = a_1 m^2/r^2 +   a_2 m^2   \exp(-c_3 r/2) . \]
So we have shown
\begin{equation}
\Pr ( N(G_m,r) >  n_m(r) )
\le a_1 c_2 \exp(- c_3r/2) .
\label{NGm}
\end{equation}
This holds for fixed $r$, but because $N(G_m,r)$ and $n_m(r)$ are decreasing in $r$ we have inclusion of events, for $j = 1, 2,\ldots $
\[ \{ N(G_m,r) > n_m(r-1) \mbox{ for some } j \le r \le j+1 \} 
\subseteq
\{ N(G_m,j) > n_m(j) \}
\]
Applying (\ref{NGm}) and summing over $j$,
 \[
\Pr ( N(G_m,r) > n_m(r-1)   \mbox{ for some } r > r_0)
\le \Phi(r_0)
\]
where $\Phi$ depends on the distribution of $\ell(e)$ but not on $m$, and 
\begin{equation}
\Phi(r_0) \downarrow 0 \mbox{ as } r_0 \to \infty.
\label{phi}
\end{equation}
Noting that $n_m(r)/m^2$ does not depend on $m$ and 
\[ \psi(r_0) := \int_{r_0}^\infty  n_m(r-1)/m^2 \ dr  \to 0 \mbox{ as } r_0 \to \infty \]
and 
$N(G_m,r) \le m^2$
we have, for all $r_0 > 0$, 
\[
\Pr \left( \int_0^\infty m^{-2} N(G_m,r)  \ dr > r_0 + \psi(r_0) \right) \le  \Phi(r_0)
\]
which, together with (\ref{phi}) and Proposition \ref{P:1}(ii), implies tightness of the sequence  $(m^{-2} L_{NUV}(G_m), \ m \ge 2)$.
\qed

The central point is that the argument depends only on some bound like (\ref{subseteq}),
which one expects to hold  very generally in FPP-like settings in dimension $> 1$.
For instance FPP on a  large family of connected random geometric graphs is studied in Hirsch, Neuh\"{a}user, Gloaguen and Schmidt \cite{hirsch}
and it seems plausible that results from that topic can be used to prove that $L_{NUV}$ is $O(n)$ on such $n$-vertex graphs.

The next example is  infinite dimensional, and the bound (\ref{Acn}) below will be the analog of the bound (\ref{subseteq}) above.

\subsection{The mean-field model of distance}
\label{sec:M-F}
Take the complete graph on $n$ vertices and assign to edges i.i.d. random weights with Exponential (mean $n$) lengths.
This ``mean-field model of distance" $G_n$ turns out to be surprisingly tractable, because 
 the smallest edge-lengths 
 $0 < \ell_1 < \ell_2 < \ldots$ 
 at a given vertex are distributed (in the $n \to \infty$ limit) 
as the points of a rate-$1$ Poisson point process on $(0,\infty)$, and as regards short edges the graph  
 is locally tree-like.
A now classical result of Frieze \cite{friezeMST} proves that the length 
$L_{MST}^{(n)}$ of the MST in this model satisfies $\Ex L_{MST}^{(n)} \sim \zeta(3) n$.
A later remarkable result of W\"{a}stlund \cite{wastlund}, formalizing ideas of  M\'{e}zard - Parisi \cite{mezard},
shows that the expected length of the TSP path in this model is asymptotically $c n$ for an explicit constant $c = 2.04....$.
Might it be possible to get a similar explicit result for the NUV length?
Corollary \ref{C:MF} below gives the correct order of magnitude by essentially the same method as above for Corollary \ref{C:grid}.
Table \ref{table:5} gives some simulation results.

\begin{table}[h!]
\centering
\begin{tabular}{rrrcc}
$n$ & $\Ex L_n$&$n^{-1} \Ex L_n$& s.d.($L_n$) & $n^{-1/2} $ s.d.($L_n$) \\
100 & 209 & 2.09 & 22 & 2.2 \\
400 &  865 & 2.14 & 41 & 2.1 \\
900 &  1954 & 2.17 & 57 & 1.9
\end{tabular}
\caption{Simulation data for lengths $L_n$ in the mean-field model.}
\label{table:5}
\end{table}

As in the previous models we expect limits of the form
\[ c := \lim_n n^{-1} \Ex L_n , \quad \sigma := \lim_n n^{-1/2} \mathrm{s.d.}(L_n) \]
and Table \ref{table:5} is loosely consistent with that.

\setlength{\unitlength}{0.7in}
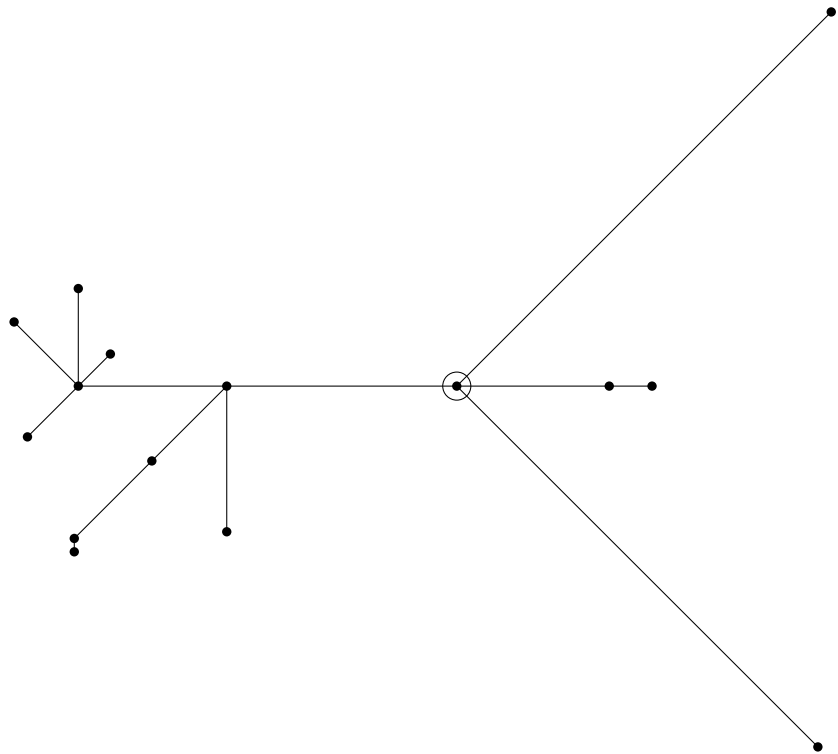
\begin{figure}
\begin{picture}(8,6)(-4,-3)
\put(0,0){\circle{0.21}}
\put(0,0){\circle*{0.07}}
\put(1.14,0){\circle*{0.07}}
\put(1.46,0){\circle*{0.07}}
\put(-1.72,0){\circle*{0.07}}
\put(-2.83,0){\circle*{0.07}}
\put(-2.59,0.24){\circle*{0.07}}
\put(-2.83,0.73){\circle*{0.07}}
\put(-3.31,0.48){\circle*{0.07}}
\put(-3.21,-0.38){\circle*{0.07}}
\put(-2.83,0){\line(1,1){0.24}}
\put(-2.83,0){\line(0,1){0.73}}
\put(-2.83,0){\line(-1,1){0.48}}
\put(-2.83,0){\line(-1,-1){0.38}}
\put(-2.28,-0.56){\circle*{0.07}}
\put(-2.86,-1.14){\circle*{0.07}}
\put(-2.86,-1.24){\circle*{0.07}}
\put(-1.72,0){\line(-1,-1){0.56}}
\put(-2.28,-0.56){\line(-1,-1){0.58}}
\put(-2.86,-1.14){\line(0,-1){0.10}}
\put(-1.72,-1.09){\circle*{0.07}}
\put(-1.72,0){\line(0,-1){1.09}}
\put(2.8,2.8){\circle*{0.07}}
\put(2.7,-2.7){\circle*{0.07}}
\put(0,0){\line(1,0){1.46}}
\put(0,0){\line(-1,0){2.83}}
\put(0,0){\line(1,1){2.8}}
\put(0,0){\line(1,-1){2.7}}
\end{picture}
\caption{Mean-field model: vertices and edges within a  ball of radius $4$ in a realization, illustrating the local tree-like property. Edges to vertices outside the ball not shown.}
\label{fig:MF1}
\end{figure}

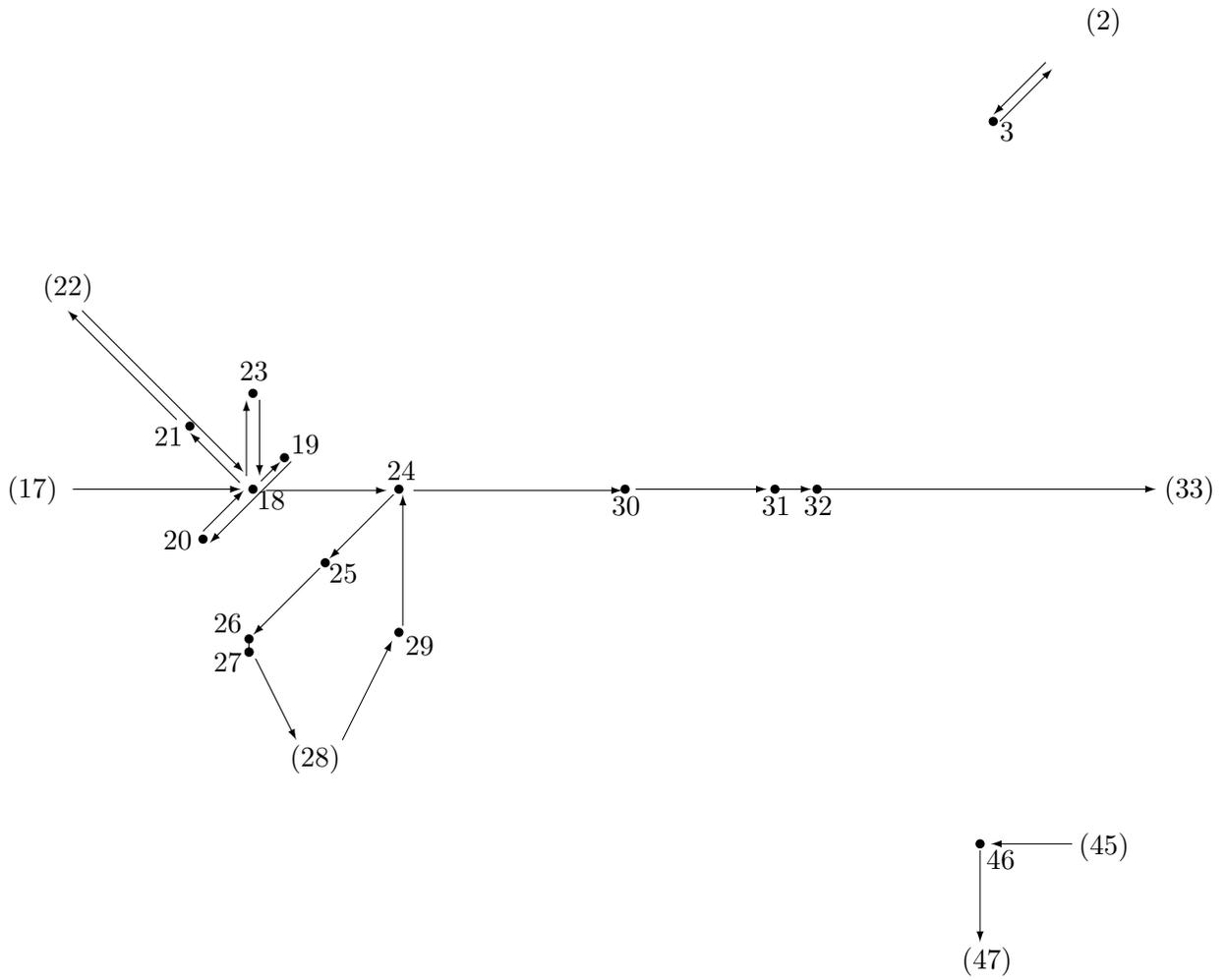
\begin{figure}
\begin{picture}(8,6)(-4,-4)
\put(2.8,2.8){\circle*{0.07}}
\put(2.85,2.65){3}
\put(3.5,3.5){(2)}
\put(2.85,2.8){\vector(1,1){0.4}}
\put(3.2,3.25){\vector(-1,-1){0.4}}

\put(2.7,-2.7){\circle*{0.07}}
\put(2.75,-2.89){46}
\put(3.4,-2.7){\vector(-1,0){0.62}}
\put(3.45,-2.78){(45)}
\put(2.7,-2.75){\vector(0,-1){0.7}}
\put(2.56,-3.65){(47)}

\put(0.08,0){\vector(1,0){1.0}}
\put(0,0){\circle*{0.07}}
\put(-0.1,-0.2){30}
\put(1.14,0){\circle*{0.07}}
\put(1.04,-0.2){31}
\put(1.12,0){\vector(1,0){0.30}}
\put(1.46,0){\circle*{0.07}}
\put(1.36,-0.2){32}
\put(1.44,0){\vector(1,0){2.6}}
\put(4.1,-0.06){(33)}

\put(-1.72,0){\circle*{0.07}}
\put(-1.61,-0.01){\vector(1,0){1.59}}
\put(-1.81,0.07){24}
\put(-1.72,-1.09){\circle*{0.07}}
\put(-1.67,-1.26){29}
\put(-1.69,-1.04){\vector(0,1){1.00}}

\put(-1.76,-0.04){\vector(-1,-1){0.49}}
\put(-2.32,-0.60){\vector(-1,-1){0.51}}
\put(-2.28,-0.56){\circle*{0.07}}
\put(-2.25,-0.71){25}
\put(-2.86,-1.14){\circle*{0.07}}
\put(-3.13,-1.09){26}
\put(-2.86,-1.24){\circle*{0.07}}
\put(-3.13,-1.39){27}
\put(-2.81,-1.29){\vector(1,-2){0.31}}
\put(-2.15,-1.91){\vector(1,2){0.38}}
\put(-2.55,-2.11){(28)}

\put(-2.83,0){\circle*{0.07}}
\put(-2.80,-0.15){18}
\put(-2.73,-0.01){\vector(1,0){0.92}}

\put(-4.2,0){\vector(1,0){1.28}}
\put(-4.7,-0.06){(17)}

\put(-2.59,0.24){\circle*{0.07}}
\put(-2.54,0.27){19}
\put(-2.77,0.06){\vector(1,1){0.15}}
\put(-2.54,0.21){\vector(-1,-1){0.62}}

\put(-3.21,-0.38){\circle*{0.07}}
\put(-3.51,-0.46){20}
\put(-3.21,-0.32){\vector(1,1){0.31}}
\put(-2.93,0.05){\vector(-1,1){0.38}}

\put(-3.31,0.48){\circle*{0.07}}
\put(-3.58,0.33){21}
\put(-3.41,0.53){\vector(-1,1){0.83}}
\put(-4.43,1.48){(22)}
\put(-4.13,1.36){\vector(1,-1){1.23}}

\put(-2.88,0.1){\vector(0,1){0.58}}
\put(-2.78,0.68){\vector(0,-1){0.58}}

\put(-2.83,0.73){\circle*{0.07}}
\put(-2.93,0.83){23}

\put(-2.86,-1.24){\circle*{0.07}}
\put(-2.86,-1.14){\line(0,-1){0.10}}
\end{picture}
\caption{Mean-field model: in the Figure \ref{fig:MF1} realization, the  NUV walk within the ball and entrance-exit edges.  
Vertices numbered according to order in an NUV walk started outside the ball, with vertices outside the ball in parentheses.}
\label{fig:MF2}
\end{figure}

\newpage

As in section \ref{sec:square}, by considering the length as $L_n(G_n,V_n)$  for a uniform random starting vertex $V_n$,
we can consider the variance decomposition
\[
\var L_n = \var \Ex(L_n \vert G_n) + \Ex \var(L_n \vert G_n) 
\]
where the first term represents the variability due to the random graph
and  the second term represents the variability due to the starting vertex. 
In simulations with $n = 100$ the former variance term is around 30 times larger than the second term,
consistent with the general conjectures (section \ref{sec:3levels}) that the initial state $v$ typically has little influence on $L_{NUV}(G,v)$.

We now prove the $O(n)$ upper bound in this model.

\begin{Corollary}
\label{C:MF}
For the mean-field model of distance $G_n$, the sequence  $(n^{-1} L_{NUV}(G_n), \ n \ge 2)$ is tight.
\end{Corollary}
To prove this, we first record a simple estimate.
\begin{Lemma}
\label{L:Hs}
Let $Z_p$ have Geometric($p$) distribution.
Let $Z^*_p$ coincide with $Z_p - 1$ outside an event $A$.
Let $H$ be a random subset of $[n] = \{1,2,\ldots,n\}$ distributed uniformly on size $Z^*_p$ subsets of $[n]$.
Then
\[ \Pr(A^c \mbox{ and } H \cap [s] = \emptyset ) \le  \frac{p}{1 - e^{-s/n}} .
\]
\end{Lemma}
\proof
It is standard (by comparing sampling with and without replacement) that
\[ \Pr(H \cap [s] = \emptyset \vert Z^*_p = i) \le \exp(-si/n) . \]
So
\begin{eqnarray*}
\Pr(A^c \mbox{ and }  H \cap [s] = \emptyset )& \le&  \sum_{i \ge 0} p (1-p)^i \exp(-si/n)\\
&=& \frac{p}{1 - (1-p)e^{-s/n}}\\ 
&\le & \frac{p}{1 - e^{-s/n}}.
\end{eqnarray*}
\qed

As before, for a vertex $v \in [n] = \{1,2,\ldots,n\}$ write
$B_n(v,r) = \{v^\prime : d(v,v^\prime) \le r \}$ for the ball of radius $r$ in $G_n$.
Conceptually we want to consider balls around $s$ randomly chosen vertices, but by symmetry this
is equivalent to using the first $s$ vertices, which is notationally simpler.
So define the vertex-set
\[ C_n(s,r) = \mbox{complement of } \cup_{i \le s} B(i,r) \]
and then by appending to $[s]$ every vertex in $C_n(s,r) $,
\begin{equation}
 N(G_n,r) \le s + |C_n(s,r)| , \ 1 \le s \le n .
 \label{NGn}
 \end{equation}
 Recall (see e.g. Pinsky and Karlin \cite{karlin} section 6.1.3) the {\em standard Yule process} $(Y(r), 0 \le r < \infty)$ for which $Y(r)$ has exactly Geometric($e^{-r}$) distribution.
The $n \to \infty$ limit distribution of the process 
$( | B_n(v,r)| , 0 \le r < \infty)$ over a fixed $r$-interval is well known to be this standard Yule process
(This is part of the theory  in Aldous and Steele \cite{PWIT} surrounding the
PWIT\footnote{Poisson Weighted Infinite Tree.}.)
Choosing $r_1 = \frac{1}{3} \log n$ so that $\exp(r_1) = n^{1/3}$ it is 
not difficult to use the natural coupling of the two processes to quantify this convergence to show 
\begin{quote}
the distribution of 
$( | B_n(v,r)| , 0 \le r \le r_1)$ agrees with the distribution of $(Y(r), 0 \le r \le r_1)$ outside an event $A_n(v)$  of probability $\delta_n = O(n^{-1/4}) \to 0$ 
as $n \to \infty$.
\end{quote}
For a vertex $v \in [s+1,n]$, and for $r \le r_1$,
\begin{eqnarray}
 \Pr(A^c_n(v) \mbox{ and }  v \in C_n(s,r)) &=& \Pr( A^c_n(v) \mbox{ and } B_n(v,r) \cap [s] = \emptyset)\nonumber \\
&\le& \frac{e^{-r}}{1 - e^{-s/(n-1)}} \label{Acn}
\end{eqnarray}
the inequality from Lemma \ref{L:Hs} applied to $[n] \setminus \{v\}$.
Apply this with
\[ s = s_n(r) := - (n-1) \log (1 - e^{-r/2}) \]
which is the solution of $e^{-r/2} = 1 - e^{-s/(n-1)}$, so 
\[  \Pr(A^c_n(v) \mbox{ and }  v \in C_n(s_n(r),r)) \le   e^{-r/2} . \]
Summing over $v$, from (\ref{NGn}) we can write, for  $r \le r_1$, 
\[ N(G_n,r) \le s_n(r) + X_n + Y_n(r)
\mbox{ where $\Ex X_n \le n \delta_n$ and
$\Ex Y_n(r) \le n e^{-r/2}$}.
\]
Applying Markov's inequality separately to the two terms on the right side of the first inequality above,
\[
\Pr( N(G_n,r) > s_n(r) + n \delta^{1/2}_n + n e^{-r/4})  \le \delta^{1/2}_n + e^{-r/4} , \ r \le r_1 .
\]
As in the proof of Corollary \ref{C:grid}
we can use monotonicity to convert this fixed-$r$ bound to a uniform bound over a ``medium" interval $r_0  \le r \le r_1$:
\[
\Pr( N(G_n,r) > s_n(r-1) + n \delta^{1/2}_n + n e^{-(r-1)/4}
\mbox{ for some } r_0 \le r \le \lfloor r_1 \rfloor 
)  \le  \delta^{1/2}_n \log n + 5 e^{-r_0/4} . 
\]
Because $s_n(r) \approx n e^{-r/2} $ over the interval of interest,
\[ n^{-1} \int_{r_0}^{r_1} (s_n(r-1) + n \delta^{1/2}_n + n e^{-(r-1)/4})  \ dr
\le K e^{-r_0/4}  + \delta_n^{1/2} \log n \]
for some constant $K$,
and so
\[ \Pr \left( n^{-1} \int_{r_0}^{r_1} N(G_n,r) \ dr > Ke^{-r_0/4}  + \delta_n^{1/2} \log n \right)
\le   \delta^{1/2}_n \log n + 5 e^{-r_0/4} . 
\]
For the tail of the integral,  the diameter $\Delta$ of $G_n$ is known (Janson \cite{janson123})  to be asymptotically $3 \log n$ 
and so by monotonicity of $N(r)$
\[n^{-1}  \int_{r_1}^{\Delta} N(G_n,r)  \ dr
= O( n^{-1} \cdot N(G_n,r_1)  \cdot \log n) 
\to 0
\mbox{ in probability}. 
\]
We will show below that
\begin{equation}
\Ex N(G_n,r_1) = O(n^{11/12}) .
\label{Nshow}
\end{equation}
Because $ \delta_n^{1/2} \log n \to 0$ and $n^{-1}   N(G_n,r) \le 1$ for $r \le r_0$, 
these bounds establish tightness of  the sequence
\[ n^{-1}  \int_{0}^{\Delta/2} N(G_n,r)  \ dr, \ \ n \ge 2 \]
which by Proposition \ref{P:1}(ii) implies 
 the sequence  $(n^{-1} L_{NUV}(G_n), \ n \ge 2)$ is tight.  
 
 To outline a proof of (\ref{Nshow}), take expectation in (\ref{NGn}) to get
 \begin{equation}
 \Ex N(G_n,r_1) \le s +  n \Pr(v \in  C_n(s,r_1)) ,    \ 1 \le s \le n 
 \label{NGn2}
 \end{equation}
 for a vertex $v \in [s+1,n]$.
 We will use this with $s = n^{3/4}$.
 Conditional on $|B_n(v,r_1)| = \beta$ we have, in order of magnitude,
 \[ \Pr(v \in  C_n(s,r_1)) \asymp (1 - \beta/n)^s \asymp \exp(- \beta s/n) .\]
 Now the distribution of $\beta$ is asymptotically Exponential with mean $e^{r_1} = n^{1/3}$,
 so by integrating over $\beta$ the unconditional probability becomes
 \[ \Pr(v \in  C_n(s,r_1)) \asymp \frac{n^{-1/3}}{n^{-1/3} + s/n} \asymp n^{-1/12} . \]
 Combining with (\ref{NGn2}) gives (\ref{Nshow}).

 \section{Final Remarks}
 \label{sec:remarks}

\paragraph{Analogy with the MST.}
As an algorithm, the NUV walk is somewhat similar to the greedy (Prim's) algorithm for the MST  (minimum spanning tree), in that both grow a connected graph one edge at at a time.
Recall that for the MST there is an intrinsic criterion for whether a given edge $e$ is in the MST 
\begin{quote}
$e$ is in the MST if and only if there is no alternative path between the endpoints of $e$, all of whose edges are shorter than $\ell(e)$.
\end{quote} 
This enables a martingale proof (Kesten and Lee \cite{kestenMST}) of the central limit theorem for the length $L_{MST}$ within the Euclidean model (complete graph on random points in the square)
which we will discuss in section \ref{sec:square}.
There is no such intrinsic criterion for the NUV walk, so
to improve the order-of-magnitude result (Corollary \ref{C:2} below) for $L_{NUV}$ in that model
one would need some other kind of control over the geometry of the set of points visited before each step.
Also, as noted in section \ref{sec:M-F}, in the 
``mean-field model of distance" the exact asymptotic constants for the lengths of the TSP tour and the MST are known: can they also be calculated for the NUV walk?

\paragraph{Local weak convergence.}
 Our results are conceptually merely consequences of Proposition  \ref{P:1},
 and further progress would require some other technique.
 One possible general approach is via local weak convergence 
 (Aldous and Steele \cite{PWIT}, 
 Benjamini and Schramm \cite{B-S}).
 Our three specific models each have local weak convergence  limits 
 (complete graph on  a Poisson point process on the infinite plane with Euclidean distance;
 i.i.d. edge-lengths on the infinite lattice; the PWIT)
 and intuitively the conjectured limits 
 $\lim_n n^{-1} \Ex L_n$ 
are the mean step-lengths in an appropriately defined NUV walk on the limit infinite graph.
Can this intuition be made rigorous?

In fact one expects the limits in our models to be  {\em collections} of disjoint doubly-infinite walks which cover the infinite graph.
 This relates to a longstanding folklore problem:
for the NUV walk on the complete-graph Poisson point process on the infinite plane, estimate the number of never-visited vertices in the radius-$r$ ball, as $r \to \infty$. 
See Bordenave,  Foss and  Last \cite{bordenave} for discussion.

\paragraph{Restrictions on local behavior of paths.}
For another possible direction of analysis, consider the Figure \ref{Fig:1} sketch of one possible trajectory for the NUV path through a given ball.
In general there will be many possible trajectories, depending on the graph outside the ball, but can one find restrictions on the possibilities,
extending the obvious restriction:
\begin{quote}
 if two vertices are each other's nearest neighbor, then every NUV walk, after visiting the first, immediately visits the second.
 \end{quote}
 Intuitively, for $1 \ll r_1 \ll r_2$, given the subgraph in the ball $B(v^*,r_2)$, in a random graph there will typically be only a few possibilities for the NUV trajectory within $B(v^*,r_1)$.

\paragraph{Variance of $L_{NUV}$?}
A final issue involves the variance of $L_{NUV}$ in random graph models.
We expect  order $n$ ``each other's nearest neighbor" pairs, and then the randomness of edge-lengths suggests that the contribution 
to variance of $L_{NUV}$ from these edges alone must be at least order $n$ (in our conventional scaling). 
However our small-scale simulation results in Tables \ref{table:7} and \ref{table:5} cast some doubt on this conjectured lower bound.

 \bigskip
 \paragraph{Acknowledgements.} I thank  three anonymous referees for helpful comments.

  \bigskip
 \paragraph{Competing interests.}  The author declares none.


\begin{thebibliography}{10}

\bibitem{PWIT}
Aldous, D.J. and Steele, J.M. (2004).
\newblock The objective method: probabilistic combinatorial optimization and
  local weak convergence.
\newblock In {\em Probability on discrete structures}, volume 110 of {\em
Encyclopaedia of Mathematical Sciences}, pages 1--72. Springer, Berlin.

\bibitem{me-FPP}
Aldous, D.J. (2016).
\newblock Weak concentration for first passage percolation times on graphs and
  general increasing set-valued processes.
\newblock {\em ALEA. Latin American Journal of Probability and Mathematical Statistics}, 13(2):925--940.

\bibitem{me-ES}
Aldous, D.J. (2021).
\newblock {Exploring Endless Space}.
\newblock In preparation.

\bibitem{auff}
Auffinger, A., Damron, M., \& and Hanson, J. (2017).
\newblock {\em 50 years of first-passage percolation}, volume~68 of {\em
  University Lecture Series}.
\newblock American Mathematical Society, Providence, RI.

\bibitem{B-S}
Benjamini, I. \& Schramm, O. (2001).
\newblock Recurrence of distributional limits of finite planar graphs.
\newblock {\em Electronic Journal of Probabability}, 6:no. 23, 13.

\bibitem{bordenave}
Bordenave C., Foss, S., \& Last, G. (2011).
\newblock On the greedy walk problem.
\newblock {\em Queueing Systems} 68:333--338.

\bibitem{cover}
Ding, J., Lee, J.R., \& Peres, Y. (2012).
\newblock Cover times, blanket times, and majorizing measures.
\newblock {\em Annals of Mathematics (2)}, 175(3):1409--1471.

\bibitem{friezeMST}
Frieze, A.M. (1985).
\newblock On the value of a random minimum spanning tree problem.
\newblock {\em Discrete Applied Mathematics}, 10(1):47--56.

\bibitem{hirsch}
Hirsch, C., Neuh\"{a}user, D.,  Gloaguen, C., \& Schmidt, V. (2015).
\newblock First passage percolation on random geometric graphs and an
  application to shortest-path trees.
\newblock {\em Advances in Applied Probability}, 47(2):328--354.

\bibitem{hougardy}
Hougardy, S. \& Wilde, M. (2015).
\newblock On the nearest neighbor rule for the metric traveling salesman
  problem.
\newblock {\em Discrete Applied Mathematics}, 195:101--103.

\bibitem{hurkens}
Hurkens, C.A.J. \&  Woeginger, G.J. (2004).
\newblock On the nearest neighbor rule for the traveling salesman problem.
\newblock {\em Operations Research Letters}, 32(1):1--4.

\bibitem{janson123}
Janson, S. (1999).
\newblock One, two and three times {$\log n/n$} for paths in a complete graph
  with random weights.
\newblock {\em Combinatorics, Probabability and Computing}, 8(4):347--361.


\bibitem{johnson}
Johnson, D.S.  \& Papadimitriou, C.H. (1985).
\newblock Performance guarantees for heuristics.
\newblock In {\em The traveling salesman problem}, Wiley-Interscience Series in
  Discrete Mathematics, pages 145--180. Wiley, Chichester.

\bibitem{kesten}
Kesten, H. (1986).
\newblock Aspects of first passage percolation.
\newblock In {\em \'{E}cole d'\'{e}t\'{e} de probabilit\'{e}s de
  {S}aint-{F}lour, {XIV}---1984}, volume 1180 of {\em Lecture Notes in Mathematics},
  pages 125--264. Springer, Berlin.

\bibitem{kestenMST}
Kesten, H. \&  and Lee, S. (1996).
\newblock The central limit theorem for weighted minimal spanning trees on
  random points.
\newblock {\em Annals of  Applied Probability}, 6(2):495--527.

\bibitem{megow}
Megow, N., Mehlhorn, K., \& and Schweitzer, P. (2012).
\newblock Online graph exploration: new results on old and new algorithms.
\newblock {\em Theoretical Computer Science}, 463:62--72.

\bibitem{mezard}
M\'{e}zard, M. \&  Parisi, G. (1986).
\newblock A replica analysis of the travelling salesman problem.
\newblock {\em Journal de Physique}, 47:1285--1296.

\bibitem{karlin}
Pinsky, M.A. \& and Karlin, S. (2011).
\newblock {\em An introduction to stochastic modeling}.
\newblock Elsevier/Academic Press.

\bibitem{rosen}
Rosenkrantz, D.J., Stearns, R.E., \& Lewis II, P.M. (1977).
\newblock An analysis of several heuristics for the traveling salesman problem.
\newblock {\em SIAM Journal on Computing}, 6(3):563--581.

\bibitem{steele}
Steele, J.M. (1989).
\newblock Cost of sequential connection for points in space.
\newblock {\em Operations Research Letters}, 8(3):137--142.

\bibitem{steelebook}
Steele, J.M. (1997).
\newblock {\em Probability theory and combinatorial optimization}, volume~69 of
  {\em CBMS-NSF Regional Conference Series in Applied Mathematics}.
\newblock Society for Industrial and Applied Mathematics (SIAM), Philadelphia,
  PA.

\bibitem{wastlund}
W\"{a}stlund, J. (2010).
\newblock The mean field traveling salesman and related problems.
\newblock {\em Acta Mathematica}, 204(1):91--150.

\bibitem{yukich}
Yukich, J.E. (1998).
\newblock {\em Probability theory of classical {E}uclidean optimization
  problems}, volume 1675 of {\em Lecture Notes in Mathematics}.
\newblock Springer-Verlag, Berlin.

\end{thebibliography}
\end{document}